\documentclass[a4paper,11pt]{amsart}
\usepackage{amscd}
\usepackage{amssymb}

\newcommand{\X}{X_{n,d}}
\newcommand{\E}{E_{n,d}}
\newcommand{\HE}{A^{\ast}(E)}

\newcommand{\DC}{{\mathcal D}}

\newcommand{\CC}{{\mathbf C}}
\newcommand{\PP}{{\mathbf P}}
\newcommand{\QQ}{{\mathbf Q}}
\newcommand{\ZZ}{{\mathbf Z}}

\newcommand{\oh}{{\mathcal{O}}}

\newtheorem{thm}{Theorem}

\newtheorem{lem}{Lemma}

\newcommand{\lav}{<\!}
\newcommand{\rav}{\!>}

\def\fra#1#2{\genfrac{}{}{0pt}{}{#1}{#2}}
 
\def\feyn#1#2#3#4{\displaystyle\Bigl(
\fra{\displaystyle\hfill #1}{\displaystyle\hfill #2}{\rangle\hskip-3pt}
\frac{\ \ }{\ \ }{\hskip-3pt\langle}
\fra{\displaystyle #4\hfill}{\displaystyle #3\hfill}\Bigr)}

\begin{document}
\title[Quantum cohomology of a Pfaffian Calabi-Yau variety]
{Quantum cohomology of a Pfaffian Calabi-Yau variety:\\ verifying 
mirror symmetry predictions}
\author[Erik N. Tj\o tta]{Erik N. Tj\o tta\\University of Bergen}
\address{Institute of Mathematics\\ University of Bergen\\ 5007 Bergen,
 Norway}
\email{erik.tjotta@mi.uib.no}
\dedicatory{To my parents on their 65th and 70th birthday.}

\begin{abstract}
We formulate a generalization of Givental-Kim's quantum
hyperplane principle. This is applied to compute the quantum
cohomology of a Calabi-Yau $3$-fold defined as the rank $4$ locus
of a general skew-symmetric $7\times 7$ matrix with  
coeffisients in $\PP^6$. The 
computation verifies the mirror symmetry predictions of
R\o dland \cite{R}.
\end{abstract}
\maketitle

\setcounter{section}{-1}

\section{Introduction}
\setcounter{page}{-1}
\pagenumbering{arabic}

The rank $4$ degeneracy locus of a general
skew-symmetric $7\times7$-matrix with
 $\Gamma(\oh_{\PP^6}(1))$-coefficients defines a non-complete
 intersection Calabi-Yau $3$-fold
$M^3$ with $h^{1,1}=1$.
We recall some results of R\o dland \cite{R} on the
 mirror symmetry of $M^3$: a potential mirror family 
$W_q$ is constructed as (a resolution of) the orbifold
$M^3_{q}/\ZZ_7$, where $M^3_{q}$ is a one-parameter family 
of invariants of a natural $\ZZ_7$-action on the space of
all skew-symmetric $7\times7$-matrices. It is shown that the
Hodge-diamond
of $W_q$ mirrors
the one of $M^3$. Further, at a  point of maximal unipotent 
monodromy\footnotemark, 
\footnotetext {There are two points with
maximal unipotent monodromy. Remarkably, the Picard-Fuchs equation
at the other point is the one found in
\cite{MR1619529} for the mirror of the 
complete intersection Calabi-Yau 3-fold
 in $G(2,7)$.}
the Picard-Fuchs operator for the periods is computed
to be (with $D=q\frac{d}{dq}$):
\begin{equation}
\label{Pic-Fuch1}
\begin{split}
(1-289q-57q^2+q^3)(1-3q)^2&D^4\\
+4q(3q-1)(143+57q-87q^2+3q^3)&D^3\\
+2q(-212-473q+725q^2-435q^3+27q^4)&D^2\\
+2q(-69-481q+159q^2
              -171q^3+18q^4)&D\\
+q(-17-202q-8q^2
             -54q^3+9q^4)&\,.
\end{split}
\end{equation}
\noindent Mirror symmetry conjectures 
that this operator is equivalent to
the operator
\begin{equation}
\label{qd}
D^2\frac{1}{K}D^2 \quad{\rm where}\quad
 K(q)=14+\sum_{d\geq 1} n_{d}d^3\frac{q^d}{1-q^d}\,,
\end{equation}
and $n_d$, the {\em instanton number of degree $d$ rational curves 
on $M^3$}, is defined \cite{MR97d:14077} using Gromov-Witten invariants by 
\[\lav p,p,p\rav^{M^3}_{d}=\sum_{k|d}k^{3}n_{k}.\]  
\noindent We shall prove the conjecture.

\begin{thm}
\label{one}
The differential operators (\ref{Pic-Fuch1}) and (\ref{qd}) are 
equivalent under mirror transformations. That is:

Let $I_0,I_1,I_2,I_3$ be a basis of solutions to (\ref{Pic-Fuch1})
with holomorphic solution $I_0=1+\sum_{d\geq 1}a_dq^d$ and logarithmic
solution $I_1=\operatorname{ln}(q)I_0+\sum_{d\geq 1}b_dq^d$. Then
\[\frac{I_0}{I_0},\frac{I_1}{I_0},\frac{I_2}{I_0},\frac{I_3}{I_0},\]
is a basis of solutions for (\ref{qd}) after change of coordinates
$q=\operatorname{exp}(\frac{I_1}{I_0})$.
\end{thm}

Our approach follows closely the work of Givental 
\cite{MR97e:14015,MR1653024} 
for complete intersections in toric manifolds, and Batyrev, 
Ciocan-Fontanine,
Kim, van Straten \cite{MR1619529,BCKV} for complete 
intersections in partial flag manifolds. It builds on the 
following three observations:
\begin{enumerate}
\item[i)]
A well-known construction identifies the degeneracy locus
$M^3$ with the vanishing locus of a section of a 
vector bundle on a Grassmannian manifold (see Section~2).
It is crucial, for us,
 that this vector bundle
decomposes into a direct sum of vector bundles $E\oplus H$, 
where $H$ is again a direct sum of line bundles. 
\item[ii)]The quantum hyperplane principle of B. Kim \cite{Kim} extends
to relate the $E$-restricted quantum cohomology with the 
$E\oplus H$-restricted one. This is formulated as a general principle
in Section~1.
\item[iii)]The $E$-restricted quantum
cohomology can be effectively computed using localization techniques and
WDVV-relations. An application of the quantum hyperplane principle 
then yields
Theorem~1. The computations are carried out in Section~2.
\end{enumerate}

\noindent {\bf Acknowledgement.} I would like to thank S.-A. Str\o mme
for bringing observation i) to my attention. Thanks also to G.~Jamet,
 B.~Kim, E.~R\o dland and D.~van Straten for helpful discussions.
Thanks to J.~Stienstra for spotting a crucial misprint in an earlier
version. 
An
extra thank you to E.~R\o dland again for letting me use his 
{\texttt Maple} scripts.\footnotemark 
\footnotetext{{\texttt http:/{/}www.math.uio.no/\~{}einara/Maple/Maple.html.}} 
I am grateful for the hospitality of 
l{'}Institut de Mathematiques de Jussieu where this work was done.
The author
has been supported by a grant from the Research Council of Norway.
 
\section{Gromov-Witten theory}
 We begin by recalling some 
basic results on $g=0$ Gromov-Witten invariants before stating the 
quantum hyperplane principle. Our approach is the algebraic 
one following \cite{MR95i:14049}. We refer the reader to 
 \cite{MR98m:14025,MR1677117} 
for a fuller account and references.

\smallskip

\noindent {\bf Frobenius rings.} Let $X$ be a smooth projective
variety
over $\CC$.
 Unless otherwise specified, we only consider even 
dimensional cohomology with rational coefficients. In fact we will
work with a further restriction: if $E$ is a vector bundle on $X$,
and $Y$ is the zero-set of a regular section of $E$, then we are
mainly interested in the cohomology classes on $Y$ that are pulled back
from $X$. These are represented by the graded Frobenius ring $\HE$ 
with
 \begin{equation}
\label{quo}
A^p(E)\colon=H^{2p}(X,\QQ)/\operatorname{ann}(E_0)
\end{equation}
and non-degenerate pairing
$\lav\gamma_1,\gamma_2\rav^{E}\colon=\int_X
\widetilde{\gamma}_1\widetilde{\gamma}_2E_0$, where $E_0$
is the top Chern class of $E$, and $\widetilde{\gamma}_i$ 
 denotes a lift of $\gamma_i$
to $A^{\ast}(X)$.

Let $A_{1}(E,\ZZ)$ be the dual of $A^{1}(E,\ZZ)/torsion$. We will
 identify 
$A_{1}(E,\ZZ)$ with the image of the natural inclusion 
\begin{equation}
\label{curveclass}
A_{1}(E,\ZZ)\to A_{1}(X,\ZZ).
\end{equation}

\smallskip

\noindent {\bf Moduli space of stable maps 
\cite{MR97d:14077,MR98m:14025}.}
  Let $(C,s_{1},\ldots,s_{n})$ 
be an algebraic curve of arithmetic genus $0$
 with at worst nodal singularities 
and $n$ nonsingular marked points. A map $f\colon C\to X$ is stable 
if all contracted components are stable (i.e. each irreducible 
component contains at least three {\em special} points, where 
special means marked or singular). For  $d\in A_{1}(X,\ZZ)$, let 
$\X$ denote the coarse moduli space (or Deligne-Mumford stack) of stable 
maps with $f_{\ast}[C]=d$. 
If $X$ is {\em convex}, that is $H^1(C,f^{\ast}TX)=0$ for all stable maps,
 then  $X_{n,d}$
is an orbifold (only quotient singularities) of complex dimension 
\begin{equation}
\label{exp.dim}
\dim_{\CC} X+\int_{d} c_{1}(X) +n-3\,.
\end{equation}

Of great importance to the theory are some natural maps on
the moduli space of stable maps.
For $i=1,\ldots,n$, let
$e_{i}\colon \X\to X$ be the map obtained by evaluating 
stable maps at $s_{i}$, and let
$\pi_{i}\colon X_{n,d}\to X_{n-1,d}$ be the map which 
forgets the marked point $s_i$. Also of significance 
are certain gluing maps
which stratify the boundaries of the moduli spaces.
In the stack theoretic framework, the diagram
\begin{equation}
        \begin{CD}
                {X_{n+1,d}} @>{e_{n+1}}>>X \\
                @V{\pi_{n+1}}VV \\
                   {X_{n,d}}
        \end{CD}
\end{equation}
along with sections $s_{i}\colon X_{n,d}\to X_{n+1,d}$ 
defined by requiring $e_i=e_{n+1}\circ s_i$, is identical to the
{\em universal} stable map.

\smallskip

\noindent {\bf Gromov-Witten invariants 
\cite{{MR93e:32028,MR96m:58033,MR95i:14049,MR99d:14011,MR98e:14022,MR98i:14015}}.}  
Suppose $E$ is a convex vector bundle (i.e $H^1(C,f^{\ast}E)=0$ 
for all
stable maps).  Base change theorems in \cite{MR57:3116} imply that 
${\pi_{n+1}}_{\ast}e_{n+1}^{\ast}E$ is a vector bundle on $\X$ 
with fibers $H^0(C,f^{\ast}E)$. 
 Let $\E$ denote 
the top Chern class of
 ${\pi_{n+1}}_{\ast}e_{n+1}^{\ast}E$ and let $c_{i}$ denote the first Chern
class of the line bundle $s_{i}^{\ast}\omega_{\pi_{n+1}}$ on $X_{n,d}$, 
where $\omega_{\pi_{n+1}}$ is the relative sheaf of differentials.
 Let $c$ be an indeterminate. 

A system of {\em $E$-restricted Gromov-Witten invariants} for $X$ is the
 family of multilinear functions $\lav\,\, \rav^{E}_{d}$ on 
$A^{*}(E)[[c]]^{\otimes n}$, defined for all $n\geq 0$ and 
$d\in A_{1}(E,\ZZ)$ by 
\begin{equation}
\label{gw}
\lav P_1\gamma_1,\ldots,P_n\gamma_n\rav_{d}^{E}\colon=
\int_{[\X]}\prod_{i=1}^{n}P_i(c_i)e_{i}^{\ast}
(\widetilde{\gamma}_i)\E\,,
\end{equation}
where  $\gamma_i\in\HE$, $P_i\in\QQ[[c]]$, and $[\X]$ is the 
{\em virtual fundamental class} of
dimension (\ref{exp.dim}). 

There is an
 exact sequence of vector bundles
\begin{equation}
\label{ker}
0\to ker\to {\pi_{n+1}}_{\ast}e_{n+1}^{\ast}E\to e^{\ast}_iE\to 0,
\end{equation}
where the right-hand map is obtained by evaluating sections at the
i-th
marked point.
This implies that $E_{n,d}$ is divisible by $E_0$ in $A^{\ast}(X_{n,d})$,  
hence the
invariants (\ref{gw}) are independent of the chosen
lifts $\widetilde{\gamma}_i$.

When $X$ is convex, then $[\X]$ is simply the 
fundamental class of $\X$. If $Y\subset X$ is cut out by 
a regular section of $E$, then

\begin{equation}
\label{class}
  j_{\ast}\sum_{i_{\ast} d'=d}[Y_{n,d'}]=\E \cdot [\X]\,,
\end{equation}

\noindent where the map $j\colon Y_{n,d'}\to \X$ is induced from 
the inclusion map
$i\colon Y\to X$.
 
 Let 
$\{\Delta_i\}$, $\{\Delta^{i}\}$ denote a pair of homogeneous
 bases 
of $\HE$ such that
$\lav\Delta_i,\Delta^{j}\rav^{E}=\delta_i^{j}$, and let
 $T_i\in\HE[[c]]$. The natural maps on the moduli space of stable
 maps respect the virtual classes, hence
 induce important relations on GW-invariants. Among
 these are:

\smallskip

\noindent {\em Divisor equation.} For $p\in A^{1}(E)$ we have
\begin{equation*}
\begin{split}
\lav p,T_1,\ldots,T_n\rav_{d}^{E}=&\\
(\int_{d}\widetilde{p})\lav T_1,\ldots,T_n\rav_{d}^{E}+&
\sum_{i=1}^{n}\lav T_1,\ldots,pT_i/c,\ldots,T_n\rav_{d}^{E}.
\end{split}
\end{equation*}
{\em WDVV-relation.} Denote\footnotemark
\footnotetext{We use the Einstein summation convention.}
\[\feyn {T_1}{T_2}{T_3}{T_4}_d \colon=\sum_{d_1+d_2=d}
\lav T_1,T_2,\Delta_i\rav_{d_1}^{E}\lav \Delta^{i},T_3,T_4\rav_{d_2}^{E}.
 \]
\noindent Then, 
\[\feyn {T_1}{T_2}{T_3}{T_4}_d=\feyn {T_1}{T_3}{T_2}{T_4}_d.\]
{\em Topological recursion relation (TRR).}
\[\lav T_1,T_2,T_3\rav_{d}^{E}=\sum_{d_1+d_2=d}
\lav T_{1}/c,\Delta_i\rav_{d_1}^{E}\lav
\Delta^{i},T_2,T_3\rav_{d_2}^{E}.\]

The above equations are subject to some restrictions:
 If $d=0$ we must assume
 that $n\geq 3$ in the
divisor equation. Further, all undefined correlators appearing
above are set to
$0$, except for
$\gamma_1,\gamma_2\in \HE$ we define 
\[\lav \gamma_1/c,\gamma_2\rav_{0}^{E}\colon=
\lav \gamma_1,\gamma_2\rav^{E}.\]

Let
$\{p_i\}$ be a {\em nef} (i.e. pairs non-negatively with all effective
curve
classes in $A_{1}(E,\ZZ)$) basis for $A^{1}(E,\ZZ)/torsion$, 
and let $\{q_i\}$
 be formal 
homogeneous parameters such that
$\sum_{i}\operatorname{deg}(q_{i})p_{i}=c_{1}(X)-c_{1}(E)$
modulo $\operatorname{ann}(E_0)$. 
The WDVV-relations imply the associativity of the {\em quantum
  product} defined by
\[\Delta_i\ast_{E}\Delta_j\colon=
\sum_{d,k}\lav \Delta_i,\Delta_j,\Delta_k\rav_d^{E}q^d\Delta^k,\]
where $q^d=\prod_i q_{i}^{\int_{d} \widetilde{p}_{i}}$. Note that 
the product is homogeneous
with the chosen grading.

\smallskip

\noindent {\bf Quantum hyperplane principle.}  
Let $\hbar$ be a formal homogeneous variable of degree $1$.
Let $e^{E}_{1\ast}$ be the map 
induced from the push-forward $e_{1\ast}$ by passing to the quotient
(\ref{quo}). Use 
$E^{'}_{1,d}$ to denote the top Chern class of the kernel in (\ref{ker}),
thus $E_{1,d}=E_0E^{'}_{1,d}$. 
Consider the following degree 0 vector in $\HE[[q,\hbar^{-1}]]$:
\begin{equation}
\label{sol}
\begin{split}
J_{E}\colon&=
e^{p\operatorname{ln}(q)/\hbar}
\sum_d q^{d}e^{E}_{1\ast}(\frac{E^{'}_{1,d}}{\hbar(\hbar-c_1)})\\
&=e^{p\operatorname{ln}(q)/\hbar}\sum_d q^{d}\lav
\frac{\Delta_{i}}{\hbar(\hbar-c)}\rav_{d}^{E}\Delta^{i},
\end{split}
\end{equation}
\noindent where 
$p\operatorname{ln}q=\sum_i \operatorname{ln}(q_i) p_i$,
 and the convention 
$\lav\frac{\Delta_{i}}
{\hbar(\hbar-c)}\rav_{0}^{E}\Delta^{i}=1$ is used.
Suppose $H=\oplus L_i$ is a sum of 
 convex line bundles on $X$.
If $c_1(X)-c_1(E\oplus H)$ is nef, the quantum hyperplane principle
suggests an explicit relationship between
 $J_E$ and $J_{E\oplus H}$ via the following
adjunct in $A^{\ast}(E\oplus H)[[q,\hbar^{-1}]]$:
\begin{equation}
\label{adjunct}
I^{H}_{E}\colon=e^{p\operatorname{ln}(q)/\hbar}\sum_d q^{d}H_{d}
e^{E\oplus H}_{1\ast}(\frac{E^{'}_{1,d}}{\hbar(\hbar-c_1)})
\end{equation}
where \[H_d\colon = 
\prod_i \prod_{m=1}^{\int_d c_{1}(L_{i})}(c_{1}(L_{i}) +
m\hbar),\]
\noindent and the $q_i$'s are regraded so that
$\sum_{i} \operatorname{deg}(q_{i})p_{i}=c_{1}(X)-c_{1}(E\oplus H)$
 modulo $\operatorname{ann(E_0)}$. Hence, 
by assumption, all $\operatorname{deg}(q_i)\geq 0$. 
The precise statement is:

\smallskip

\noindent {\em The vectors $I^{H}_{E}$ and $J_{E\oplus H}$ coincide up to
 a {\em mirror transformation} of the following type:
\begin{enumerate}
\item[i)]multiplication by $a\operatorname{exp}(b/\hbar)$ where $a$ and
$b$ are homogeneous $q$-series of degree $0$ and $1$ respectively,
\item[ii)]coordinate changes 
\[\operatorname{ln}q_i\mapsto\operatorname{ln}q_i + f_i\,,\]
where $f_i$ are homogeneous $q$-series of degree $0$ without constant term.
\end{enumerate}
Further, the mirror transformation is uniquely determined by the
first two coefficients in the $\hbar^{-1}$-Taylor expansion of 
$I^{H}_{E}$ and $J_{E\oplus H}$}.

\begin{thm}
\label{qlp}
If $X$ is a homogenous space and $E$ is
equivariant with respect to a maximal torus action on $X$, then
the quantum hyperplane conjecture as formulated above is true.
\end{thm}
\begin{proof}
The proof in \cite{Kim} for $\operatorname{rank}(E)=0$ extends
with minor modifications to the general case.
\end{proof}

\noindent {\em Remark~1.}
An early version of this principle appeared in
\cite{MR96g:32037}.
 Givental formulated and proved the
 $\operatorname{rank}(E)=0$
 case for toric manifolds \cite{MR97e:14015,MR1653024}. 
 The above formulation  when
$\operatorname{rank}(E)=0$
is due to B.~Kim \cite{Kim}.  Extending the conjecture of B.~Kim
we expect that the principle holds for more general $X$.
In \cite{T} the conjecture was tested on a non-convex 
non-toric manifold.

\smallskip

\noindent {\em Remark~2.} An analogue generalization of the
hyperplane principle in \cite{MR1672116,Kim2}  for 
{\em concavex} $H$ 
can be formulated with $E$ convex/concave. 
The proof in \cite{Kim}
extends to cover these cases when $X$ is 
homogenous. See also \cite{MR99e:14062}.

\medskip

\noindent {\bf Differential equations.}
The vector $J_E$ encodes all the $E$-restricted one-point GW-invariants. 
Reconstruction using TRR \cite{MR1685628} shows that these are 
determined by two-point GW-invariants without
$c$'s. 
This is organized nicely in terms of differential equations
\cite{MR97d:58038,MR97e:14015}.
 Consider the {\em quantum differential
equation}
\begin{equation}
\label{qdiff}
\hbar q_{k}\frac{d}{dq_{k}}T
=
p_{k}\ast_{E}T, \quad k=1,\ldots,\operatorname{rank}(A^1(E)),
\end{equation}
where $T$ is a series in the variables $\operatorname{ln}q_i$ 
and $\hbar^{-1}$ with
coefficients from $A^{\ast}(E)$.  The WDVV-relations imply that 
the system is solvable. An application of the divisor equation and
 TRR \cite{T,MR1685628} 
 shows that
$J_{E}=\lav S,1\rav^{E}$ for fundamental solution $S$ of (\ref{qdiff}).
In particular, if $c_1(X)-c_1(E)$ is 
{\em positive} the hyperplane principle, when true,
yields an algebraic representation of the ({\em quantum}) $\DC$-module
 generated by 
$J_{E\oplus H}$.

\smallskip

\noindent {\em Remark~3.} A useful application of Theorem~2 is the
following:
for partial flag manifolds $F=F(n_1,\ldots,n_r,n)$ 
with universal sub-bundles $U_{n_i}$ and quotient bundles 
 $Q_{n_i}$, one may consider vector bundles $E$ that are direct sums
of bundles of type
\[\wedge^{p_1}U^{\vee}_{n_{i_1}}\otimes \wedge^{p_2}Q_{n_{i_2}}\otimes
S^{p_3}U^{\vee}_{n_{i_3}}\otimes S^{p_4}Q_{n_{i_4}}.\]
\noindent If $c_1(F)-c_1(E)$ is positive, the $E$-restricted quantum
cohomology
can in principle be computed using localization techniques, and the 
theorem will yield the quantum $D$-module for the nef
(and in particular the Calabi-Yau) cases of type 
$E\oplus H$, with $H$ decomposable.

\section{The Pfaffian variety}
Let $V$ be a vector space of dimension $7$ and consider the
projective space $P=\PP(\wedge^{2}V)$ with 
 universal $7\times 7$ skew-symmetric
linear map 
\[\alpha\colon V_{P}^{\vee}(-1)\to V_{P},\]
where $V_P$ denotes the trivial vector bundle on $P$ with fiber $V$. 
Define $M\subset P$ as the locus where 
$\operatorname{rank}\alpha \leq 4$.
The scheme structure is determined by the Pfaffians of the diagonal
$6\times 6$-minors. The variety is locally Gorenstein of codimension
$3$ in $P$ with canonical sheaf $\oh_{M}(-14)$. Its singular locus, 
which is the rank $2$ degeneracy locus of $\alpha$, is of codimension $7$ 
in $M$ \cite{MR56:11983}. This implies that the intersection 
$M^k=M\cap\PP^{k+3}$ with a
general linear sub-space $\PP^{k+3}$ in $P$ is of dimension $k$,
has canonical sheaf $\oh_{M^k}(3-k)$, and is smooth when $k\leq 6$. 

We recall a classic construction for degeneracy loci
(see for instance \cite{MR99d:14003}, Example 14.4.11).
 Let $G=\operatorname{Grass}_4(V)$ be the Grassmannian of 
$4$-planes in $V$ 
with universal exact sequence  
\[0\to U\to V_G \to Q\to 0\,.\]
 Pulling everything back to $PG=P\times G$, we regard $\wedge^2U(1)$ as a sub-bundle
of $\wedge^2V_{PG}(1)$, where the twists are with
respect to $\oh_P(1)$. 
The map $\alpha$ induces
a regular section $\overline{\alpha}$ of the convex rank $15$ quotient bundle
on $PG$
\[A\colon=\wedge^2V_{PG}(1) /\wedge^2U(1).\]

\begin{lem}
\label{SAS}
The zero-scheme $V(\overline{\alpha})\subset PG$ projects 
birationally onto $M$. 
Moreover, the projection is isomorphic over the non-singular locus of $M$.
\end{lem}

For a general linear sub-space $\PP^{k+3}$ in $P$, let 
$(A^k, \overline{\alpha}_k)$ denote the pull-back of the pair 
$(A, \overline{\alpha})$ to $\PP^{k+3}\times G$. By Lemma~\ref{SAS}, 
$V(\overline{\alpha}_k)$ projects isomorphically to $M^k$ for $k\leq 6$.
 
We are now
set
to compute the quantum $\DC$-module of the Calabi-Yau variety $M^3$
using Theorem~\ref{qlp}. This can in principle be done from
any of the $A^k$-restricted $(k\geq 4)$ GW-theories.
We provide details for the case $E=A^{6}$.

First, we need to determine the cohomology ring $A^{\ast}(E)$.
 Pull-backs to $\PP^9\times G$ of the Chern
classes $p=c_1(\oh_P(1))$ and $\gamma_i=c_i(Q), i=1,2,3$,  generate the
$\QQ$-algebra $A^{\ast}(\PP^9\times G)$, 
and induce generators of $A^{\ast}(E)$.
 
\begin{lem}
\label{cohom}
\begin{enumerate}\item[]
\item[i)] The $\QQ$-algebra $A^{\ast}(E)$ is generated by $p$ and $\gamma_2$ with
top degree monomial values as follows:
\begin{center}
\begin{tabular}{cccc}
$p^6=14$&$p^4\gamma_2=28$&$p^2\gamma_2^2=59$&$\gamma_2^3=117$
\end{tabular}
\end{center}
\noindent In particular the Betti numbers of $A^{\ast}(E)$
 are $[1,1,2,2,2,1,1]$.
\item[ii)] We have the relation $\gamma_1=2p$.
\end{enumerate}
\end{lem}

\begin{proof} Computed using {\texttt Schubert}.\footnotemark
\footnotetext{A {\texttt{MAPLE}} package for
enumerative geometry written by S.~Katz and S.-A.~Str{\o}mme.
 Software and documentation available at
{\texttt http:/{/}www.math.okstate.edu/\~{}katz/schubert.html}.}
\end{proof}

Our choice of basis $\{\Delta_i\}$ for $A^{\ast}(E)$ is the following:
\begin{equation*}
\{1,p,p^2,\gamma_2,p^3,p\gamma_2,p^4,p^2\gamma_2,p^5,p^6\}.
\end{equation*}
Rather than to work directly with the solutions (\ref{sol}) and
 (\ref{adjunct}) we prefer to work with
their governing differential equations.
The quantum differential equation (\ref{qdiff})
is determined by the two-point numbers $\lav
 \Delta_i,\Delta_j\rav_{d}^{E}$
for all $d$ in $A_1(E,\ZZ)\simeq\ZZ$. A simple dimension 
count shows that
 these are $0$
 unless 
\begin{equation}
\label{codim}
\operatorname{codim}(\Delta_i)+\operatorname{codim}(\Delta_j)=5+3d\,.
\end{equation} 

\begin{lem}
\label{2_point} Values of $d\geq 1$ GW-invariants satisfying (\ref{codim}) 
are as follows:\\

\centering
\begin{tabular}{ccc}
$\lav p^2,p^6\rav_{1}^{E}=238$&$\lav p\gamma_2,p^5\rav_{1}^{E}=2044$&
$\lav p^2\gamma_2,p^2\gamma_2\rav_{1}^{E}=6617$\\
$\lav \gamma_2,p^6\rav_{1}^{E}=504$&$\lav p^4,p^4\rav_{1}^{E}=1568$&\\
$\lav p^3,p^5\rav_{1}^{E}=980$&$\lav p^4,p^2\gamma_2\rav_{1}^{E}=3220$&
$\lav p^5,p^6\rav_{2}^{E}=9800$
\end{tabular}
\end{lem}
\begin{proof}
It follows from Lemma~\ref{cohom} that the map
 (\ref{curveclass}) identifies $d$ with the
curve class $(d,2d)$ in $A_1(\PP^9\times G,\ZZ)$, so the GW-invariants are
integrals over $(\PP^9\times G)_{(d,2d)}$. 

{\em Localization.} Consider the standard action of 
$T=(\CC^{*})^{10}\times (\CC^{*})^7$ 
on $\PP^9\times G$. Since the integrands are polynomials in Chern classes of 
equivariant vector bundles with respect to the
induced $T$-action on $(\PP^9\times G)_{(d,2d)}$, they may
be evalued using Bott's residue formula
(see \cite{MR96j:14039,MR97d:14077}). 
As there are only finitely many fixed points and curves in $\PP^9\times G$, 
the formulae involved are similar to the ones found in
 \cite{MR97d:14077}. Details
are left to the reader. The two-point integrals with
 $d=1$ were evaluated in this manner. 

{\em Reconstruction.} The number $\lav p^5,p^6\rav_2^{E}$
is however more easily obtained from $d=1$ numbers using the 
following WDVV-relations:
\begin{equation}
\label{ass}
\feyn {p}{p^2\gamma_2}{\gamma_2}{p^5}_2\,, \feyn {p}{p^4}{\gamma_2}{p^5}_2\,,
\feyn {p}{p^3}{\gamma_2}{p^6}_2\,, \feyn {p}{p\gamma_2}{\gamma_2}{p^6}_2\,.
\end{equation}
\noindent In fact, a further analysis\footnotemark 
\footnotetext {For instance using {\texttt Farsta}, a computer program written
  by A.~Kresch, available at 
{\texttt http:/{/}www.math.upenn.edu/\~{}kresch/computing/farsta.html}.
}
shows that the $3$-point $d=1$
numbers appearing in the equations (\ref{ass}) are in turn determined by 
the $2$-point $d=1$ numbers above.
\end{proof}
\noindent {\em Remark.} Employing a description in \cite{EPS} of the class in 
$\operatorname{Grass}_2(\wedge^2V)$
of lines on $M$, the $d=1$ GW-invariants which only involve powers
of $p$ can be computed without using Bott's formula.

\smallskip

Consider the differential equation (\ref{qdiff}), 
 with invariants as in Lemma~\ref{2_point} and $\hbar=1$. Denote $q=q_1$
 and $D=q\frac{d}{dq}$. 
By reduction
we find the (order $10$, degree $5$)-differential equation 
$P(D)=\sum_d q^{d}P_d(D)=0$, with $P_d$ as below, for $J_E(\hbar=1)$.

{\small \begin{align*}
P_0=&3D^7(D-1)^3,\\
P_1=&D^3(194D^7-776D^6+1072D^5-1405D^4-1716D^3-1272D^2-414D-51),\\
P_2=&343D^{10}-1715D^9+3185D^8-58593D^7-55484D^6-460D^5+10697D^4+\\
    &1850D^3-896D^2-480D-96,\\
P_3=&-99127D^7+22736D^5-11772D^4-34797D^3-31654D^2-13495D-2175,\\
P_4=&-19551D^4-39102D^3-31360D^2-11524D-1430,\\
P_5=&343(D+1).
\end{align*}}

Let $H=3\oh_P(1)$ and assume $\hbar =1$. The adjunct $I^{H}_{E}$ 
is obtained by
correcting the $q^d$-coefficients of $J_E$ with the class
 $H_d=\prod_{m=1}^{d}(p+m)^3$.
A reformulation of this transformation on the corresponding
differential equations takes the same form\footnotemark. That
is, the differential operator
\footnotetext{This is a general principle when 
$\operatorname{rank}A^1(E)=1$. It is easily proved 
using recursion formulas for
 solutions of differential equations \cite{MR96g:32037,T}.}
\begin{equation}
\label{diff}
\sum^5_{d=0}q^{d}P_d\prod_{m=1}^{d}(D+m)^3
\end{equation}
\noindent annihilates $I^{H}_{E}$. Recall that the commutation rule
is $Dq-qD=q$.
 If we factor out the ``trivial'' term $D^3(D-1)^3$ from the 
left of (\ref{diff}) and 
re-organize the terms we recover the Picard-Fuchs operator
 (\ref{Pic-Fuch1}).

\begin{proof}[Proof of Theorem~\ref{one}]
From (\ref{class}) it follows that
$\lav p,p,p\rav^{M^3}_{d}=\lav p,p,p\rav^{A^{3}}_{d}$.
Using (\ref{qdiff}) it is straightforward to check
that $D^2\frac{1}{K}D^2$ is the differential operator
governing $J_{A^{3}}$ (see for instance \cite{T}). The rest follows from 
Theorem~\ref{qlp}.
\end{proof}

The first five curve numbers are:\\

\centering
\begin{tabular}{ccccc}
$n_1=588$&$n_2=12103$&$n_3=583884$&$n_4=41359136$&$n_5=3609394096$
\end{tabular}

\end{document}